\documentclass[12pt]{article}
\usepackage{amsmath}
\usepackage{amssymb}
\usepackage[dvips]{epsfig,graphics} 

\begin{document}

\newenvironment{proof}{\par\noindent{\bf Proof:\/}}{\hfill$\square$\vskip .3cm}
\newtheorem{thm}{Theorem}[section]
\newtheorem{lemma}[thm]{Lemma}
\newtheorem{prop}[thm]{Proposition}
\newtheorem{obs}[thm]{Fact}
\newtheorem{defi}[thm]{Definition}
\newtheorem{cor}[thm]{Corollary}
\newtheorem{fact}[thm]{Fact}

\newcommand{\g}{\mathfrak g}
\newcommand{\s}{\mathfrak s}
\newcommand{\toro}{\mathfrak t}
\newcommand{\const}{\rm const}
\newcommand{\rank}{\rm rank~}
\newcommand{\rat}{\rm rat~}
\newcommand{\im}{\rm im~}
\newcommand{\symm}{\rm Symm}
\newcommand{\Diag}{\rm Diag}
\newcommand{\diag}{\rm diag}
\newcommand{\ve}{\varepsilon}
\newcommand{\vp}{\varphi}
\newcommand{\vt}{\vartheta}
\renewcommand{\o}{\overline}
\newcommand{\wt}{\widetilde}
\newcommand{\R}{{\mathbb R}}
\newcommand{\Z}{{\mathbb Z}}
\newcommand{\C}{{\mathbb C}}
\newcommand{\Q}{{\mathbb Q}}

\author{Andrea Giacobbe\\
Department of Mathematics\\
University of Maryland at College Park\\
College Park, MD 20742, U.S.A.}
\title{Convexity of multi-valued momentum map}
\maketitle

\begin{abstract}

Given a Hamiltonian torus action the image of the momentum map is a convex polytope; this famous convexity theorem of Atiyah, Guillemin and Sternberg still holds when the action is not required to be Hamiltonian. Our generalization of the convexity theorem states that, given a symplectic torus action, the momentum map can be defined on an appropriate covering of the manifold and has as image a convex polytope along a rational subspace times its orthogonal space. We also prove stability of this property under small perturbations of the symplectic structure. The technique developed allows us to extend the result to any compact group action and also to deduce that any symplectic torus action, with fixed points, on a manifold of double the dimension of the torus, is Hamiltonian.
\end{abstract}

\tableofcontents

\section{Introduction}

When a group $G$ acts on a symplectic manifold $M$ it is sometimes possible to define a map from $M$ to the dual Lie algebra of $G$; this map is called momentum map. The momentum map is obtained defining a linear map from the Lie algebra of $G$ to the closed 1-forms of $M$; these 1-forms are called (multi-valued) Hamiltonians. In order to define the momentum map we need the Hamiltonians to be exact 1-forms. If this condition is verified the action is called Hamiltonian. We refer to \cite{g-s} for definitions and properties of the momentum map.

The famous convexity theorem of M.F. Atiyah, V. Guillemin and S. Sternberg states that the image of the momentum map for a Hamiltonian torus action on a closed manifold is a convex polytope. Many works followed this result but, until very recently\footnote{I was informed that a result identical to Theorem \ref{thm of multi-valued convexity} was announced in a short communication by Y. Benoist \cite{bns98}.}, they all dealt with the case of Hamiltonian actions only.
  
To relax such global request is relevant, as was pointed out by S.P. Novikov \cite{nvk}. Hamiltonian dynamic can be originated by Hamiltonians that are closed 1-forms in important physical problems, as happens in presence of Dirac monopoles.

Our working environment consists of a closed connected symplectic
manifold, equipped with a symplectic action of a compact connected Lie
group. The main result in this work describes the structure of the
image of the momentum map when the group acting is a torus. The result
of Atiyah, Guillemin and Sternberg is quite unaltered. The difference
from that theorem sits on the fact that the momentum map cannot be
described as a function on the manifold $M$, because closed 1-forms
are not functions on $M$. It is although possible to define a map $\wt
\mu$ from $\wt M$, a covering of the manifold $M$, to the dual Lie algebra of $T$. The image of $\wt\mu$ is a bounded, convex polytope along a rational subspace of the Lie algebra times the orthogonal space. The technique developed allows us to deduce an interesting result. If a torus of half the dimension of $M$ acts with fixed points, then the action is Hamiltonian. 

It is very natural, at this point, to consider the work of F. Kirwan and to obtain a convexity result for any compact connected Lie group. Kirwan \cite{krw84} proved that the intersection of the image of the momentum map with a Weyl chamber in a Cartan subalgebra is a convex polytope. We can obtain the equivalent result for non Hamiltonian actions.

The last section is devoted to a stability result for Hamiltonian actions. We show that a small perturbation of the 2-form does not effect the exactness of Hamiltonian 1-forms arising from a group action. The result is interesting since there is no apparent reason that would prevent Hamiltonian vector fields to be perturbed into locally Hamiltonian vector fields.

\section{Normalization}

As anticipated in the introduction we are given a closed, connected, symplectic manifold $(M^{2N}, \Omega)$, with an action of a compact, connected Lie group $G^n$. Everywhere in this paper the action is assumed to be symplectic. 

Recall that any smooth action on a manifold yields an anti-Lie algebra map$$
\g\stackrel{X_\bullet}{\longrightarrow} Vect(M),
$$
called infinitesimal action; the vector fields $X_\xi$ are called fundamental vector fields. The symplectic 2-form allows to define a linear map to the de Rham 1-forms
$$
\g\stackrel{\tau_\bullet}{\longrightarrow} A^1_{dR}(M).
$$
The action is symplectic exactly when the image of the map $\tau_\bullet$ is contained in the closed de Rham 1-forms $Z^1_{dR}(M)$.  The 1-forms associated to a group action are called Hamiltonians. 

The convexity properties of the momentum map are based on two important facts. The first of them is that, in a neighborhood of any point, the momentum map can be explicitly written in local coordinates. 

This is proven by first, linearizing the action in a neighborhood of a point fixed by the action, and then changing coordinates so that the symplectic form is put in canonical form. The approach we adopt is that of Guillemin and Sternberg \cite{g-s82}. 

\begin{defi} If $p$ is a point fixed by the action of $G$ then the representation 
$$
G\rightarrow GL(T_pM)\simeq GL(2N,\R), \qquad g\mapsto d_pg
$$
is called isotropy representation.
\end{defi}

A proof of the three following lemmas can be found in \cite{g-s} or \cite{g-s82}.
 
\begin{lemma}
Let $(M^{2N},\Omega)$ be a symplectic manifold, $G$ a compact connected Lie group acting on $M$, $p\in M$ a point fixed by the action. There exist an equivariant open neighborhood of $p$ and coordinates on that neighborhood such that the action of $G$, in these coordinates, is the isotropy action.
\end{lemma}

\begin{lemma}\label{moser and weinstein}
Let $M^{2N}$ be a manifold, $G$ a compact connected Lie group acting on $M$, $p\in M$ a point fixed by the action, $\Omega_0$ and $\Omega_1$ two $G$-invariant symplectic 2-forms such that $\Omega_0(p)=\Omega_1(p)$. There exist an equivariant local diffeomorphism
$$
(M,p,\Omega_0)\longrightarrow (M,p,\Omega_1).
$$
\end{lemma}

In this lemma ``equivariant local diffeomorphism''  means that there exist an equivariant neighborhood of $p$ and a diffeomorphism to an equivariant neighborhood of $p$ that fixes $p$ and respects the action. This diffeomorphism also pulls back the symplectic form $\Omega_1$ to $\Omega_0$.

We now apply these results to the $n$-dimensional torus. The fact that irreducible torus representations are one dimensional over $\C$ allows us to obtain a very special choice of coordinates in a neighborhood of a point fixed by the action. Locally a torus action is exactly a family of harmonic oscillators.

\begin{lemma}\label{normalization}
Let $(M^{2N},\Omega)$ be a symplectic manifold, $T^n$ an $n$-dimensional torus acting on $M$, $p\in M$ a point fixed by the action. In a neighborhood of $p$ we can define $N$ complex valued local coordinates $z^j=x^j+iy^j$ such that:

\noindent
$\bullet$ the symplectic form, in these local coordinates, has expression
$$
\Omega = dx^1\wedge dy^1 +\cdots +dx^N\wedge dy^N.
$$ 

\noindent
$\bullet$ The action of the torus is the linear action associated to the representation
$$
\begin{array}{ccc}
T^n&\longrightarrow&\Diag\cap U(N)\\
(\Theta_1,...,\Theta_n)&\mapsto&\diag(e^{i(\alpha_1^1\Theta_1+\cdots+\alpha^n_1\Theta_n)},...,e^{i(\alpha_N^1\Theta_1+\cdots+\alpha^n_N\Theta_n)}),
\end{array}
$$
where $\alpha_j$ are real characters of the isotropy representation, and hence belong to $2 \pi \Z^n \subset \R^n$. The sign of the characters is determined by the fact that the symplectic form induces an orientation on the coordinate planes. 

\noindent
$\bullet$ The momentum map, in these local coordinates, has expression
$$
\begin{array}{ccc}
\C^N & \stackrel{\mu}{\longrightarrow} & \toro^* \\
(z^1,...,z^N) & \mapsto & \alpha_1 \frac{|z^1|^2}2+\cdots+\alpha_N
\frac{|z^N|^2}2,
\end{array}
$$
where the $\alpha_j$ are now thought as real infinitesimal weights of the isotropy representation.
\end{lemma}

It is easy to show that the Hamiltonian nature of the action yields abundance of fixed points. In our setting the action is not Hamiltonian hence we need a local description of the momentum map in a neighborhood of a generic point $p\in M$. 

Given a point $p\in M$ its stabilizer, denoted by $T_p$, is a non necessarily connected subtorus whose dimension we denote by $n_p$. For simplicity we define $m_p=n-n_p$. Unfortunately, the fact that the action is not necessarily isotropic\footnote{An action is isotropic if the fundamental vector fields span an isotropic subspace of the tangent space.} forces us to introduce an auxiliary integer number $r_p$ which satisfies $2r_p\leq m_p$. This integer turns out to be the rank of the skew-symmetric matrix of commutation relations for the Hamiltonians.

\begin{lemma}\label{normalization on non fixed points}
Let $(M^{2N},\Omega)$ and $T^n$ be as above, let also $p$ be any point in $M$. In a neighborhood of $p$ we can find complex valued local coordinates $w^1,...,w^{m_p-r_p}$ and $z^1,...,z^{N-m_p+r_p}$ such that:

\noindent
$\bullet$ the symplectic form, in these local coordinates, has expression 
$$
\Omega=\sum d\Re z^j\wedge d\Im z^j + \sum d\Re w^j\wedge d\Im w^j.
$$

\noindent
$\bullet$ The action, for $\Theta \in S^1$ and for $|\Phi|$ small is
\begin{equation}\label{action}
\begin{array}{c}\
(\Theta_1,...,\Theta_{n_p},\Phi_1,..., \Phi_{m_p})(z^1,...,z^{N-m_p+r_p},w^1,...,w^{m_p-r_p})=\\
=(e^{i(\alpha_1^1\Theta_1+\cdots + \alpha^{n_p}_1\Theta_{n_p})}z^1,..., e^{i(\alpha_{N-m_p}^1\Theta_1+\cdots+\alpha^{n_p}_{N-m_p+r_p}\Theta_{n_p})}z^{N-m_p+r_p},\\
,w^1+\Phi_1+i\Phi_2,..., w^{r_p}+\Phi_{2r_p-1}+i\Phi_{2r_p}, w^{r_p+1}+\Phi_{2r_p+1},...,w^{m_p-r_p}+\Phi_{m_p}),
\end{array}
\end{equation}
where $\alpha_j$ are real characters of the isotropy representation of $T_p$ in $\Diag(N-m_p+r_p,\C)\cap U(N-m_p+r_p)$. The sign of the characters is determined by the symplectic form.

\noindent
$\bullet$ The momentum map, in these local coordinates, has expression\begin{equation}\label{momentummap}
\begin{array}{c}
\C^{N-m_p+r_p}\times \C^{m_p-r_p} \longrightarrow\toro_p^*\times(\toro_p^*)^\perp, \\
( z^1,...,z^{N-m_p+r_p},w^1,...,w^{m_p-r_p})\mapsto(\alpha_1 \frac{|z^1|^2}2+\cdots +\alpha_{N-m_p+r_p} \frac{|z^{N-m_p+r_p}|^2}2,\\
,\Re w^1,\Im w^1,...,\Re w^{r_p},\Im w^{r_p},\Im w^{r_p+1},...,\Im w^{m_p}),
\end{array}
\end{equation}
where the $\alpha_j$ are now thought as real infinitesimal weights of
the isotropy representation of $T_p$ in $GL(T_pM)$.
\end{lemma}

In this lemma we used the fact that $\toro$ has a distinguished subgroup isomorphic to $\Z^n$: all periodic vector fields of integer period. Given a basis in the Lie algebra $\toro$ there is an Ad-invariant metric: the Euclidean metric that has the given basis as orthonormal. This metric allows an identification between the Lie algebra of $T$ and its dual. Under this identification is understood any non canonical map. In many statements we will use this observation.

The coordinates of $T^n$ we used to write the action are first order coordinates for a chosen basis of $\toro$. The coordinates $\Phi$ are not periodic coordinates of $T$, while the coordinates $\Theta$ are. Formula \ref{action} can be written for a different choice of basis. In this case the expression of the real infinitesimal weights, and consequently the image of the momentum map, changes by left multiplication by an element of $SL(n_p,\Z)\times GL(m_p,\R)$.

\begin{proof}
Let $\vt_1,...,\vt_{n_p}$ be an integral basis of $\toro_p$; complete it with integral vectors $\vp_1,...,\vp_{m_p}$ to a basis of $\toro$. In a neighborhood of $p$ the vectors $\vp_j$ generate $m_p$ independent, commuting, locally Hamiltonian vector fields $X_{\vp_j}$. The closed forms $\tau_{\vp_j}$ define a family of $m_p$ functions $x^j$. The commuting relations of such functions originate a constant skew-symmetric matrix. 

By a linear change of coordinates in $GL(m_p,\R)$ --- hence not preserving integrality --- we can replace the vectors $\vp_j$ with a new family of vectors whose associated 1-forms have commuting relations that are represented by the rank $2r$ skew symmetric matrix 
$$
\left(
\begin{array}{cc|c}
0&-1&\\
1&0&\\
\hline
&&0
\end{array}
\right).
$$

We still denote $\vp_j$ the new vectors and $x^j$ the functions associated to them. The vector fields $X_{\vp_j}$ commute, hence we can find functions $y^j$, $j=1,...,2N$ such that $\partial_{y^j}=X_{\varphi_j}$ for $j=1,...,m_p$. The $2m_p-2r$ functions $x^1,...,x^{2r},x^{2r+1},...,x^{m_p},y^{2r+1},...,y^{m_p}$ are independent and, by a simple calculation $\{y^l,x^j\}=\nabla^sx^j(y^l)=X_{\varphi_j}(y^l)=\partial_{y^j}y^l=\delta^{lj}$. Hence the commuting relations for these functions can be summarized by the matrix
$$
\left(
\begin{array}{cc|cc}
0&-1\\
1&0\\
\hline
&&0&-1\\
&&1&*\\
\end{array}
\right),
$$
which has blocks of dimension $2r$ and $2m_p-2r$.

The construction we just did gives no informations on what the matrix $*$ could be, but we can mimic a technique in \cite{ar'd74}: consider the manifold $V=\{x^j=0, y^l=0|j=1,...,m_p, l=2r+1,...,m_p\}$. $V$ has an induced symplectic structure because its defining functions have non degenerate Poisson brackets. Define the remaining $2(N-m_p+r_p)$ functions by choosing canonical (complex) coordinates in a neighborhood of $p$ in $V$, $z^1,...,z^{N-m_p+r_p}$, and extending them to a neighborhood of $V$ as constant along the flow of the symplectic gradient of the functions $x^1,...,x^{2r},x^{2r+1},...,x^{m_p},y^{2r+1},...,y^{m_p}$.\footnote{$\{\nabla x^j,\nabla y^l\}\cap TV=0$. In fact, if a vector $a_j \nabla x^j +b_l\nabla y^l$ is in $TV$ then, by bracketing it against the function $x^l$, we obtain that $b_l=0$. Using then $\{y^l,-\}$, we prove the vanishing of the coefficients $a_l$. The dimensions add up so that the flows of the symplectic vector fields, $\nabla y^l$ for $l=2r_p+1,...,m_p-r_p$ and $\nabla x^j$ for $j=1,...,m_p$, fill injectively a tubular neighborhood of $V$.}

Finally we can replace the functions $y^{2r+1},...,y^{m_p}$ with functions in involution. Note in fact that $\{y^j, y^k\}$ are functions of $x^1,...,x^m_p$ only, so we can obtain the new coordinates replacing the function $y^l$ by the function $y^l+f^j(x^1,...,x^{m_p})$ for appropriate $f^j$. This can be done in view of the fact that $\sum_{j,k}\{y^j,y^k\}dx^j \wedge dx^k$ is a closed 2-form, and hence it is locally the differential of a 1-form $\sum_j f^j dx^j$; this defines the functions $f^j(x^1,...,x^{m_p})$. Finally define the functions $w^j=x^j+ix^{r_p+j}$ for $j=1,...,r_p$ and $w^l=x^{r_p+l}+iy^l$ for $l=1,...,m_p-r_p$.    

It is easy to verify that $T_p$ acts as the identity on the coordinates $w^j$. Hence the action of $T_p$ has isotropy representation that can be reduced to a representation in $\Diag(N-m_p+r_p,\C)\cap U(N-m_p+r_p)$.
\end{proof}

\section{Convexity for torus action}

\subsection{Local convexity}

In the previous section we proved that the momentum map admits an explicit expression in a neighborhood of a critical point; our goal is to prove global convexity of the image of the momentum map also when it is not possible to globally define it as a function on $M$, because the $T$-action is not Hamiltonian. In the rest of the work we will often refer to the local momentum map at $p\in M$, with this name we design the momentum map $\mu_p$ defined in a coordinate chart $U\subset M$, $p\in U$, such that $\mu_p(p)=0$. The local momentum map is well defined; in fact the Hamiltonians are closed 1-forms and a coordinate chart is simply connected.

We state here an useful fact the proof of which is almost a tautology.

\begin{fact}
The rank of the local momentum map at a point $p$ ($\rank d_p\mu$) is the dimension of the orbit through $p$ ($\dim Tp$) which is also the codimension of the stabilizer of $p$ in $T$ ($\dim T -\dim T_p$).
\end{fact}
   
From formula \ref{momentummap} on page \pageref{momentummap} we can
easily deduce that the image of the local momentum map at $p$ maps the neighborhood of $p$ to the vertex at the origin of a convex wedge $W_p$. To be more precise let $p\in M$ be a critical point for the local momentum map (which is the same as assuming that its stabilizer $T_p$ is a non discrete subgroup) then the image of the local momentum map at $p$ is the convex wedge
$$
\R^{\geq}\{\alpha_j|j=1,...N-m_p+r_p\}=\{\sum
\lambda_j\alpha_j|\lambda_j\in\R^{\geq}\}\subset\toro_p^*
$$
times the vector space $(\toro^*_p)^\perp$, where $\alpha_j$ are the real infinitesimal weights found in Lemma \ref{normalization on non fixed points}. The choice of $\toro_p^*$ is non canonical but follows from a choice of an algebraic complement to $\toro_p$. 

\vskip .5cm
\centerline{\mbox{\epsfig{file=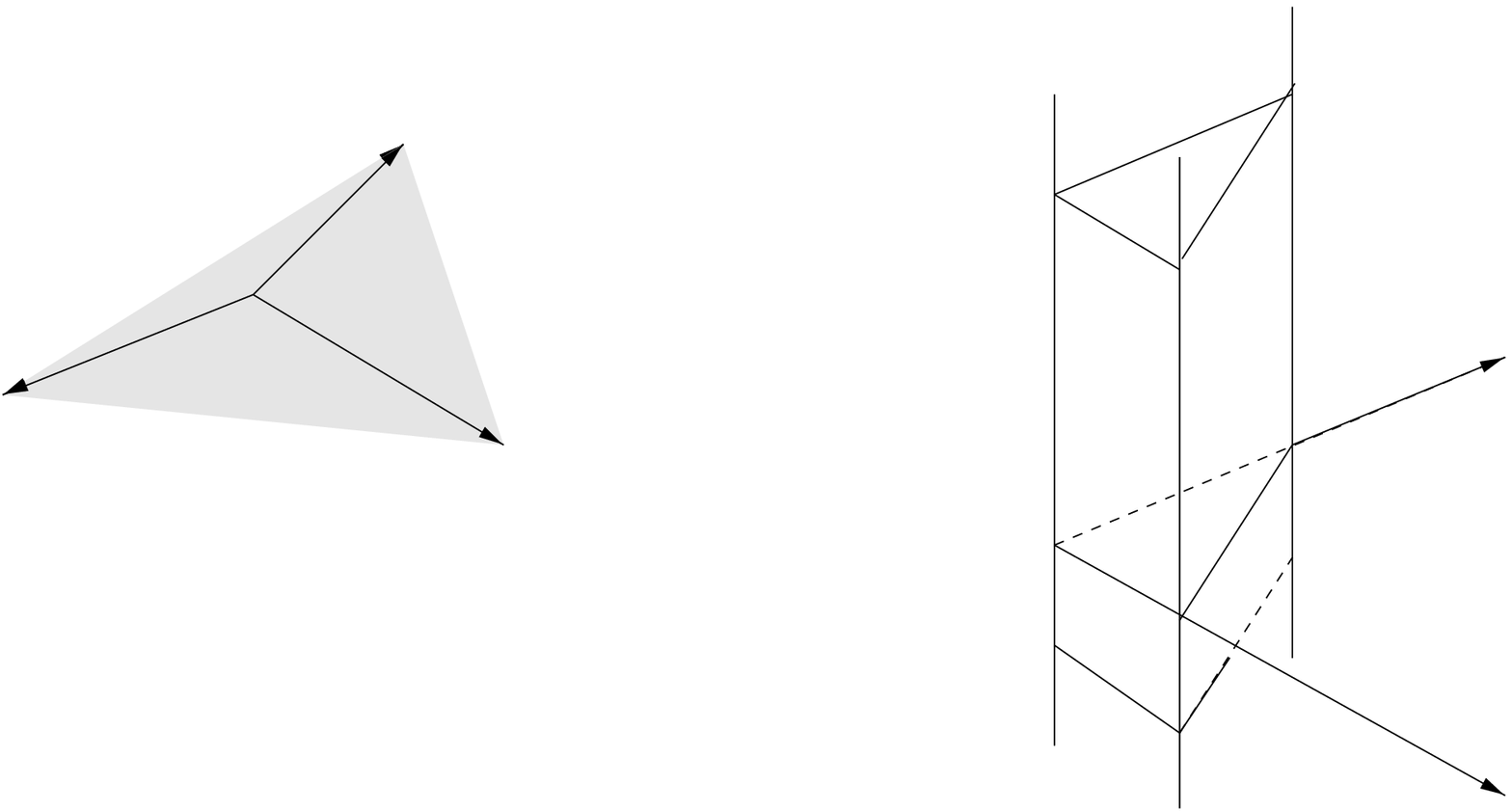,height=100pt}}}
\centerline{Fig. 1}

\begin{defi}
An action is called effective if the only element of the group that fixes all points of $M$ is the identity.
\end{defi}

\begin{prop}\label{effective action}
A torus action with fixed points is effective if and only if, at a fixed point, the real infinitesimal weights associated to the isotropy representation are generators of $\Z^n$ in the dual Lie algebra $\toro^*$.
\end{prop}

In this statement we are again using an Ad-invariant scalar product to identify covectors to vectors. To be more precise we choose an integral basis of $\toro$ that, moreover, generates $\Z^n\subset\toro$. We can hence identify elements of $\toro^*$ with elements of $\toro$. It is hence defined a distinguished lattice $\Z^n$ in $\toro^*$.

\begin{proof}
From Lemmas \ref{normalization} and \ref{normalization on non fixed points} it is obvious that the $\alpha_j$ belong to $\Z^n$. Let $\alpha_j$, for $j=1,...,N$, be the list of real infinitesimal weights at $p$. 

Assume that the $\alpha_j$ are a set of $\Z$-generators of $\Z^n$. An element of $\Theta\in T$ stabilizes all points of $M$ if and only if $\langle \alpha_j,\Theta \rangle \equiv 0$ for all $j$ ($\equiv$ means $\bmod 2\pi$) if and only if $\langle \lambda^j\alpha_j,\Theta \rangle \equiv 0$ for all $\lambda^j \in \Z$, if and only if $\Theta$ is an element of $2\pi \Z^n$, if and only if $\Theta$ is the identity of the torus. 

For the converse, if we assume that $\{\alpha_j\}_\Z$ is a strict subgroup of $\Z^n$, by eventually an $SL(n,\Z)$ change of basis we can assume such subgroup is $q_1\Z \oplus \cdots \oplus q_n\Z$ with at least $q_n$ not 1. The element of the torus $\bar \Theta=(0, \cdots, \frac1{q_n})$ stabilizes all the points of the manifold\footnote{We are using a well known a fact about group actions: generic points of $M$ have principal type (minimal stabilizer) \cite{adn}. In case of effective torus action this means that the points stabilized only by the identity of $T$ form a dense open subset of $M$; hence every open set has to contain such points.} without being 0. This contradicts effectiveness. 
\end{proof}

The previous statement, in case of an $n$-dimensional torus action on a $2n$-dimensional manifold, implies that we can choose a basis of the torus so that the real infinitesimal weights are the dual basis of an integral basis of $\toro$. In this case the image of the local momentum map, at a fixed point, is the standard tetrahedron, as was observed by T. Delzant.

\begin{defi}
Let $p$ be a point of $M$, the rank of $p$ is the codimension of its stabilizer in the torus, $\dim T-\dim T_p$. The exceptionality of $p$ is the finite group $\frac{T_p}{(T_p)_0}$.
\end{defi}

A point $p\in M$ is mapped to the boundary of wedge associated to some (hence any) local momentum map if and only if:

\noindent
$\bullet$ $p$ is stabilized by some, non discrete subgroup $T_p$.

\noindent
$\bullet$ The $\alpha_j$ arising from the isotropy representation are such that $\R^{\geq}\{\alpha_j\}$ (the convex wedge consisting of the positive linear combinations of the real infinitesimal weights $\alpha_j$) is not the whole $\R^{n_p}\simeq\toro_p^*$.

Observe that, if the dimension of the torus is half of the dimension of the space, or if the dimension of the critical manifold which $p$ belongs to is $2N-2n$, this second request is verified. 


\subsection{Brief digression on Morse-Bott 1-forms}\label{MB1-forms}

In next subsection we will use a property of closed Morse-Bott 1-forms. The main source on this subject is \cite{btt54}. The result we will use is

\begin{thm}\label{minimum->exact}
Let $\tau$ be a  closed Morse-Bott 1-form, without critical manifolds of index one, in a closed connected manifold $M$. If $\tau$ has local minimum it has exactly one minimum and it is exact.  
\end{thm}
\begin{proof}
The notion of local maximum and local minimum makes perfect sense for closed 1-forms. A proof of this fact would require a long digression hence we omit it. A partial result (for rational 1-forms) is used in \cite{t-w00}. A full proof can be found in \cite{primo}.
\end{proof}

\subsection{Half dimensional torus actions}

We apply the previous theorem to the Hamiltonians $\tau_\xi$ obtained
from a torus action. Since the map $\tau_\bullet$ is linear there is a distinguished subspace in $\toro$ which is the preimage of the exact 1-forms. Its elements are vectors associated to one-valued Hamiltonians. We denote this subalgebra by $\toro_e$.

The elements of $\toro$ not in $\toro_e$ are associated to closed non exact 1-forms that we call multi-valued Hamiltonians. The results in \cite{aty82} and \cite{g-s82} are stated in the case in which all Hamiltonians are one-valued, this condition is expressed by calling the action Hamiltonian.

The forms $\tau_\vt$, for $\vt\in\toro$, are all closed Morse-Bott 1-forms; their critical manifolds are symplectic and so even dimensional, the index of $\tau_\vt$ at a critical point is twice the cardinality of the set $\{j|\langle\alpha_j,\vt\rangle<0\}$ hence is even; it is obvious that the coindex is even too. All these statements follow from the local expression of the momentum map.

\begin{thm}\label{N=n implies Hamiltonian}
Let $M$ be a $2n$-dimensional, closed, connected symplectic manifold, $T$ be an $n$-dimensional torus acting on $M$; assume also that the action has fixed points. Then the action is Hamiltonian.
\end{thm}
\begin{proof} Let $p$ be a point fixed by the action. The real infinitesimal weights of the isotropy representation at $p$, $\{\alpha_j|j=1,...,n\}$, must $\Z$-generate $\Z^n$ and are exactly in number of $n$. By an $SL(n,\Z)$ change of coordinates in $\toro^*$ we can assume that $\alpha_j=e_j^*$, where $e_j$ is an integral basis for $\toro$ and the $*$ indicates the dual element ($e_i^*(e_j)=\delta_{ij}$). So the image of the local momentum map is the vertex at the origin of the standard tetrahedron.

Choose now any $\vt\in\toro$ that is in the first octant of $\toro$, in other words choose a vector all whose coefficients with respect to the chosen basis are positive. The form $\tau_\vt$ has a local minimum at $p$ so, it must be exact. 
\vskip .5cm
\centerline{\mbox{\epsfig{file=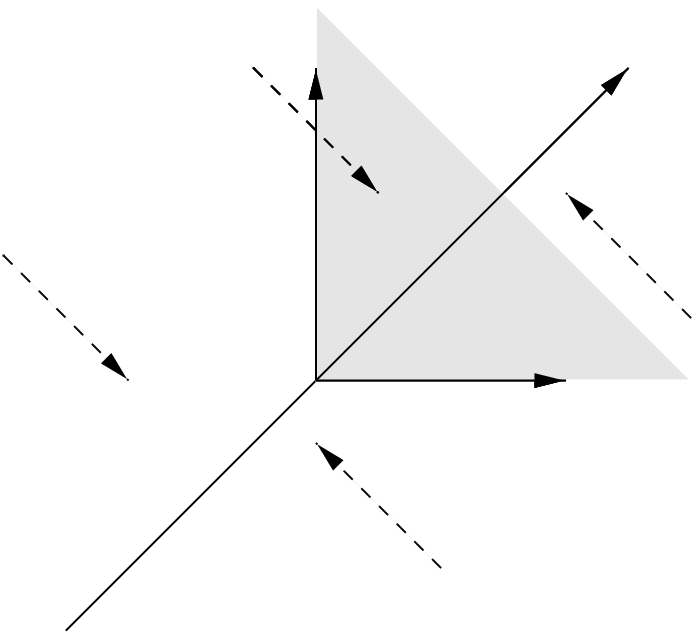,height=150pt}}}
\centerline{Fig. 2}

The subset of $\toro$ associated to exact Hamiltonians is a subvector space and all the first quadrant is mapped to exact 1-forms. Hence the action is Hamiltonian.
\end{proof}

\subsection{Rationality degree}

\begin{defi}
A closed 1-form $\tau$ defines a map 
$$
\int_\bullet\tau:H_1(M,\Z)\longrightarrow \R, \qquad \gamma \mapsto \int_\gamma\tau.
$$
The rationality degree of $\tau$ (rat $\tau$) is the dimension of the $\Q$-vector space that the image of $\int_\bullet\tau$ spans.
\end{defi}
 
Having to deal with a subspace of closed 1-forms, and not with a unique 1-form, it is not obvious how to define the rationality degree. 

\begin{prop}\label{collective kernel}
For almost any choice of $\vt$ in $\toro$ the multi-valued Hamiltonian $\tau_\vt$ induces a map from $H_1(M,\Z)$ to $\R$ that has minimal kernel. Hence, the smallest covering of $M$ in which all the pulled back Hamiltonians are one-valued is the one in which a generic Hamiltonian is one-valued.
\end{prop}

Before proving the statement we observe that the minimality condition on the kernel is more than a condition on the rationality degree, in fact 
$$
\rat\tau=\rank\left(\frac{H_1(M,\Z)}{\ker\int_\bullet\tau}\right).
$$
\begin{proof} Let $\{\vt_i|i=1,...,n\}$ be a basis of the Lie algebra $\toro$. This basis is associated to $n$ 1-forms $\tau_{\vt_i}$; to each of them is associated a finite dimensional $\Q$-vector space $V_i=\Q\{\im \int_\bullet\tau\}\subset\R$. For any $n$-tuple of real numbers $(\lambda_i)$ such that the set $\{\lambda_iV_i|i=1,...,n\}$ of $\Q$-subvector spaces of $\R$ is in direct sum, the form associated to the vector $\sum\lambda_i\vt_i$ has minimal kernel. In other words, given any such $n$-tuple, for any $\vt$ in $\toro$ 
$$
\ker \int_\bullet\lambda_i\tau_{\vt_i}\supset\ker\int_\bullet\tau_\vt.
$$
Observe that the $n$-tuples for which the family $\{\lambda_iV_i\}$ is not composed of free summands are a countable subset of a continuum space.
\end{proof}

An interesting fact, which implies an upper bound for the possible codimension of the subalgebra $\toro_e$ is the following

\begin{prop}
Assume a torus acts on a closed symplectic manifold so that there are no one-valued Hamiltonians except the one associated to the origin of $\toro$ ($\toro_e=0$). Then a generic Hamiltonian must be at least $n$-rational, where $n$ is the dimension of the torus.
\end{prop}
\begin{proof} Assume that the generic Hamiltonian has rationality degree $r$. Then we can choose a basis $\{\gamma_i|i=1,...,b_1(M)\}$ of the free part of the first homology group of $M$ so that
$$
\int_{\gamma_j}\tau_\vt=0 \quad \forall j>r, \vt\in\toro.
$$
This follows from the Theorem of Elementary Divisors and Proposition \ref{collective kernel}.

If the generic Hamiltonian was less than $n$-rational then the matrix $\int_{\gamma_i} \tau_{\vt_j}$ would have less than $n$ non zero lines. Hence there would be a relation of linear dependence among the lines. Such relation specifies a choice of an exact Hamiltonian, contradicting the hypothesis.
\end{proof}

It becomes obvious that, if an $n$ dimensional torus acts on a manifold with first Betti number $b_1(M)$ smaller than $n$, then the subtorus $\toro_e$ associated to one-valued Hamiltonians must be at least $n-b_1(M)$ dimensional.

\subsection{Global convexity}

The following definition was introduced by Atiyah in \cite{aty82}.

\begin{defi}
A function is said to be quasi-periodic if the symplectic flow it generates determines a finite dimensional torus action.
\end{defi}

\begin{prop}
On a closed symplectic manifold with symplectic torus action $\toro_e$, the subspace of vectors in $\toro$ associated to exact Hamiltonians, is a rational subspace. 
\end{prop}
With this proposition we prove that the Lie subalgebra $\toro_e$ exponentiates to a subtorus. We denote such subtorus with $T_e$.

\begin{proof} 
The thesis is equivalent to the fact that, if a quasi-periodic function generates a $T^n$-action, then $\toro_e=\toro$. The negation of this last statement can be proven to be equivalent to the following situation: a torus $T^n$, with $n>1$, acts on a closed symplectic manifold. This action is such that $\toro_e=\R\xi$ for some $\xi\in\toro$ maximally irrational (i.e. the one parametric subgroup generated by $\xi$ is dense in $T$).

Let $p$ be a point that is a maximum for the Hamiltonian $\tau_\xi$; this point must exist given that the manifold is compact and $\tau_\xi$ is exact; it is easy to prove that $p$ is also a fixed point by the action of $T$. 

The fact that $p$ is a maximum for $\tau_\xi$ implies that the real infinitesimal weights of the isotropy representation at $p$, that we denoted by $\alpha_j$, must satisfy the inequalities
$$
\langle \alpha_j , \xi \rangle \le 0 \quad \forall j=1,...,N.
$$
Some of these inequalities must be equalities, otherwise an argument as in Theorem \ref{N=n implies Hamiltonian} would imply that the action is Hamiltonian. Call $J_0$ the set of indices such that equality holds.

We plan to show that $\R^\geq\{ \alpha_j | j \in J_0 \}=\xi ^\perp$. If not there is a vector $\eta\in\xi^\perp$ (we are identifying $\toro$ with $\toro^*$) such that $\langle\alpha_j,\eta\rangle\le 0$ for all $j \in J_0$. So, for $\ve\in\R$ small,
$$
\begin{array}{llll}
\langle \alpha_j , \xi + \varepsilon \eta \rangle & < & 0 &\forall j \notin J_0,\\
\langle\alpha_j,\xi +\varepsilon\eta\rangle&\le&0&\forall j\in J_0.
\end{array}
$$
By Proposition \ref{minimum->exact} this implies that the vectors $\xi +\ve\eta$ are associated to one-valued Hamiltonians, contradicting the hypothesis that $\toro_e=\R\xi$. We hence deduce that the space $\xi^\perp$ is rational, and $\xi$ cannot be maximally irrational.
\end{proof}

\begin{thm}\label{thm of multi-valued convexity}
Let $M$ be a closed symplectic manifold, $T$ an $n$-dimensional torus acting on $M$. There is a regular covering $\widetilde M$ in which any 1-form $\tau_\vt$ pulls back to an exact form; from such covering the momentum map can be defined as a function $\wt\mu: M \rightarrow\toro^*$. The image of $\wt\mu$ is a convex polytope times a vector space. 
\end{thm}

To be more precise the set $\wt\mu(\wt M)$ is the convex polytope obtained as follows: take $T_e$, the maximal, one-valued subtorus of $T$ that acts in a Hamiltonian fashion. Its fixed manifolds $Z_j\subset M$, determine the vertices of a convex, bounded polytope $P$ in $\toro^*_e$ that is the image of the momentum map for the action of $T_e$ on $M$. We will see that the set $\wt\mu(\wt M)$ is $P\times (\toro_e^*)^\perp$ (for some scalar product).

\begin{proof} In Proposition \ref{collective kernel} we proved that the smallest covering from which it is possible to define the momentum map as a function is that in which the pull back of a generic Hamiltonian is exact. Such covering is regular, with group of deck transformations which is free Abelian of rank equal to the rationality degree of the generic Hamiltonian.

We denote the elements of the deck transformation of $\wt M/M$ by $I=(i_1,...,i_r)\in\Z^r$, the projection in $M$ of a point $\wt p\in \wt M$ will be denoted by $p$. 

\begin{lemma}\label{convex wedge}
Let $p$ be a point of $M$. The image of $\mu_p$,the local momentum map at $p$, defines a convex wedge, $W_p\subset \toro^*$ (that could well be $\toro^*$). For any point $\wt p$ in $\wt M$ above $p$ the image of the momentum map $\wt\mu$ is contained in 
$$
\wt\mu(\wt p)+W_p\subset\toro^*.
$$
\end{lemma}
\begin{proof}
The convex wedge $W_p$ is defined by some inequalities
$$
\langle\vt,-\rangle\geq 0.
$$
Each of such inequalities yields a local maximum at $p$ for the Hamiltonian $\tau_\vt$. By Theorem \ref{minimum->exact} the Hamiltonian $\tau_\vt$ is exact and has unique maximum at $p$. This implies that for any $\wt x$ in $\wt M$ and for any $\wt p$ above $p$
$$
\begin{array}{cc}
\langle\vt,\wt\mu(\wt x)\rangle &\equiv \langle\vt,\wt\mu(I\wt x)\rangle \\
&\geq \langle\vt,\wt\mu(\wt p)\rangle,
\end{array}
$$
and this shows exactly that the image of $\wt\mu$ is bounded by the set of affine equations needed to conclude.
\end{proof}

We can now prove the convexity of $\im\wt\mu$. Let $\psi_1$ and $\psi_2$ be two points in $\im\wt\mu$, assume the interval $[\psi_1,\psi_2]=\{(1-t)\psi_1+t\psi_2|t\in[0,1]\}$ is not contained in the image of $\wt\mu$. Let $\vp_\infty$ be the first point in the segment not in the interior of $\im\wt\mu$.

We can find a sequence $\wt p_n$ in $\wt M$ such that $\wt\mu(\wt p_n)=\vp_n\in [\psi_1,\vp_\infty]$ and $\vp_n$ converges to $\vp_\infty$. The sequence $p_n$ admits a subsequence $p_{n'}$ that converges to a point $p_\infty$ in $M$.

The local momentum map at $p_\infty$ maps a neighborhood $U$ of $p_\infty$ into a convex wedge $W_{p_\infty}$. The points $p_n$ belong to $U$ for all $n$ big enough. For every $n$ big enough there is an unique open set $\wt U_n\subset\wt M$ containing $\wt p_n$ and a point $\wt p_{\infty,n}$ above $p_\infty$.

The point $\vp_\infty$ must belong to all the convex wedges 
$$
\wt\mu(\wt p_{\infty,n})+W_{p_\infty},
$$
otherwise it could not be in the closure of the image of $\wt\mu$. Hence it is an attained value $\vp_infty=\wt\mu(\wt q)$. It is easy to be convinced that $W_q$ is a proper wedge of $\toro^*$ and that it does not contain the semi-line
$$
\{\vp_\infty+t(\psi_2-\vp_\infty)|t>0\}.
$$
By Lemma \ref{convex wedge} we conclude that $\psi_2$ cannot belong to the image of $\wt\mu$.
\begin{lemma}\label{3.14}
If the Hamiltonian map is such that $\toro_e=0$, then the image of the momentum map $\wt\mu$ is surjective.
\end{lemma}
\begin{proof} 
Let $\gamma_1,...,\gamma_{b_1}$ a base for the cycles of the free part of $H_1(M,\Z)$ and call 
$$
v_j=(\int_{\gamma_j} \tau_{\vt_1},..., \int_{\gamma_j} \tau_{\vt_n})^t.
$$
The $v_j$'s are $b_1$ vectors of $\R^n$. 

The hypothesis implies that such vectors are a system of generators of $\R^n$. In fact, if they were not, then the lines of $A=(v_1\cdots v_{b_1})$ would be dependent and so    
$$
\sum_i\lambda_i \int_{\gamma_j} \tau_{\vt_i}=0 \quad \forall j.
$$
This expresses a non zero, one-valued Hamiltonian. We can also think of $\lambda$ as a non zero vector orthogonal to all the $v_j$; such a vector can be found if and only if the $v_j$'s do not span $\R^n$.

Since $\wt\mu(\wt M)$ must contain all the $\Z$-subgroup generated by the vectors $v_j$; in fact $\tilde \mu(I\tilde p)=\tilde\mu(\tilde p)+\sum_{i_j}v_j$. Hence, by the just proven convexity, the image of $\wt\mu$ must be $\toro^*$.   
\end{proof}

Since the subalgebra $\toro_e$ is rational we can choose a rational algebraic complement $\toro_c$ so that not only $\toro=\toro_e\oplus\toro_c$ and $\toro^*=\toro_e^*\oplus\toro_c^*$, but also $T^n=T_e^{n_e}\oplus T_c^{n_c}$.

The projection of the momentum map $\wt\mu$ on $\toro_e^*$ factors through a function $\mu_e:M\rightarrow\toro_e^*$. The image of $\mu_e$ is a convex polytope as described in \cite{aty82} and \cite{g-s82}. This means that there are critical manifolds fixed by the $T_e$ action, $Z_j$, and $\im\mu_e={\rm co}\{\mu_e(Z_j)\}=P$. This implies that $\im\wt\mu\subset \im\mu_e\times\toro_c^*=P\times\toro_c^*$. 

We are left to show that the image of $\wt\mu$ is exactly $P\times\toro_c^*$.  Let $\wt p$ be a point of $\wt M$ that projects on $p\in M$, and let $I$ be an element of $\Z^r$, the deck transformations of $\wt M$. The element $\tilde \mu(I\wt p)\in\toro^*_e\oplus\toro_c^*$ has first component that is the same as that of $\wt\mu(\wt p)$, while the second component is translated by a vector, associated to $I$, as described in Lemma \ref{3.14}.

We can use, again, convexity of the image of $\wt\mu$ to show that $\im\wt\mu\supseteq P\times\toro_c^*$. This concludes the proof.
\end{proof}

\section{Convexity for compact group action}

\subsection{Preliminaries and notations}

In this section we prove a convexity result when the torus is replaced by any compact connected Lie group. Notions such as regular elements, Cartan subalgebra and Weyl group are used. Our notations are consistent with the ones in \cite{b-d}.
 
The most important structure result for compact Lie groups is the following. 
 
\begin{fact}\label{structure of compact group}
Let $H$ be a compact connected Lie group. There exists a finite covering of $H$ which is isomorphic to $G\times T$. $G$ is a semisimple simply connected Lie group, $T$ is a torus and is also the center up to finite extension.
\end{fact}

\begin{lemma}\label{action of the covering}
Let $H$ be a compact connected Lie group that acts on a symplectic manifold $(M,\Omega)$, let $\wt H$ be a finite covering of $H$. The momentum maps associated to the two group action on $M$ have same image.
\end{lemma}
\begin{proof}
The Lie algebras of the two groups are the same and the fundamental vector fields are the same.
\end{proof}

The action of $\wt H$ is surely not effective, but effectiveness is not going to play any role in what follows. Let $H$ be a compact connected Lie group that acts on a closed connected symplectic manifold $M$. By Proposition \ref{structure of compact group} we know that $H\simeq G\times_F T$, where $F$ is a discrete subgroup of $G\times T$. Proposition \ref{action of the covering} states that the image of the momentum map for the $H$-action is equivalent to that of the momentum map for the $G\times T$-action. Since $G$ is semisimple its action is Hamiltonian, while the action of $T$ can be multi-valued.

We need to set some more notations for what follows. In $\g$ we can choose a maximal commutative subalgebra, called Cartan subalgebra and denoted by $\mathfrak{s}$; the exponential of $\s$ is a subgroup of $G$ isomorphic to a torus, denoted by $S$. The normalizer of $S$ in $G$ induces a linear action (a finite group action) on $\mathfrak s$, this action is the well known the Weyl group action; a fundamental domain of such action is called Weyl chamber and is denoted by $\mathfrak s_+$; with $\mathfrak s_0$ we denote the interior of the Weyl chamber. 

We define as $\nu:M\rightarrow\g^*$ the momentum map for the $G$ action on $M$; $\wt\mu:\wt M\rightarrow \g^*\times\toro^*$ is the momentum map for the action of $G\times T$ on $M$.

\begin{defi}
The subset of elements of $\g$ that belong to a unique Cartan subalgebra are called regular.
\end{defi}

The set of regular elements of the Lie algebra will be denoted by $\g_{reg}$. The set $\g-\g_{reg}$ is a union of submanifolds of codimension at least 3 in $\g$. The set $\mathfrak s_0$ coincides with the set of regular elements in the Weyl chamber. 

In their work \cite{g-s82} Guillemin and Sternberg introduced the idea of ``symplectic cross-section''; the set $\nu^{-1}(\mathfrak s_0^*)=V$ is not only a submanifold of $M$ but satisfies the following properties:

\noindent
$\bullet$ $V$ is an $S\times T$-invariant symplectic submanifold of $M$.

\noindent
$\bullet$ Every orbit in $M_{reg}$ intersects $V$, where $M_{reg}$ is the open set of $M$ that $\nu$ maps to regular elements of $\g$.

\noindent
$\bullet$ The set $\nu(V)$ is dense in $\nu(M)\cap\mathfrak s_+^*$. 

The only statement not proven in \cite{g-s82} is the $T$-invariance of $V$. This follows from the fact that every Hamiltonian arising from the action of $T$ Poisson-commutes with every Hamiltonian arising from the action of $G$.

\subsection{Convexity in the Weyl chamber}

The convexity of $\im\nu\cap\mathfrak s_+^*$ was conjectured in \cite{g-s82} and was finally proven in \cite{krw84}. $\wt M$ contains a submanifold, $\wt V$, obtained as the full preimage of $V$ under the projection map. On $\wt M$ is also defined a momentum map $\wt\mu:\wt M\rightarrow \g^*\times\toro^*$. The restriction of $\wt\mu$ to $\wt V$ has image contained in $\s^*_0\times\toro^*$.

\begin{lemma}\label{4.6}
Let  $\vt\in\toro$ be such that the Hamiltonian $\tau_\vt$ is a closed, non-exact 1-form on $M$. Then the form $(\tau_\vt)_{|V}$ cannot have local maximum nor minimum in $V$.
\end{lemma}
\begin{proof}
Assume the form $(\tau_\vt)_{|V}$ has a local maximum at the point $q\in V$. We claim that $\tau_\vt$ must then have a local maximum at $q$. In fact the form $\tau_\vt$ is orthogonal to vectors of the form $X_\xi$, for any $\xi\in\g$ (in other words $\langle\tau_\vt,X_\xi\rangle=0$ for any $\xi\in\g$); hence the Hessian of the form $\tau_\vt$, at the point $q$, has same signature and same rank as the Hessian of $(\tau_\vt)_{|V}$ at the same point.
\end{proof}

\begin{prop}
Let $M$ be a closed symplectic manifold, $G\times T$ a compact
connected Lie group acting on $M$, assume that the torus $T$ acts on
$M$ so that all its infinitesimal vector fields are locally but not
globally Hamiltonian; let $\wt \mu$ and $\nu$ be the maps defined in Subsection 4.1. The set $\im\wt\mu\cap(\mathfrak{s}_+^*\times\toro^*)$ is equal to $P\times\toro^*$, where $P=\im\nu\cap\mathfrak s^*_+$ is a convex polytope contained in the Weyl chamber $\mathfrak{s}_+^*$.
\end{prop}
\begin{proof}
In \cite{krw84} was proven that the set $\nu(M)\cap\s_+^*=\o{\nu(V)}$ is a convex polytope $P$, and that $\nu(V)=\nu(M)\cap\s_0^*$. We just have to show that
\begin{equation}\label{noncommutativeconvexity}
\left(\wt \mu(\wt V)\right)={\nu(V)}\times \toro^*.
\end{equation}

We plan to use the fact that $V$ is an $S\times T$-manifold, and that the
momentum map $\wt\mu$, restricted to $\wt V$, has values in $\s^*\times\toro^*$.

Assume \ref{noncommutativeconvexity} is false, let $\vp$ be a point in $\nu(V)\times\toro^*$ and $\vt^*\in\toro^*$ such that the line $\vp+\R\vt^*$ is not all contained in $\wt\mu(\wt V)$. The line in question must intersect $\wt\mu(\wt V)$, since the projection on $\s^*$ of $\wt\mu(\wt V)$ is $\nu(V)$.

With a line of proof similar to that in Theorem \ref{thm of multi-valued convexity}\footnote{There is a difference here: $V$ is non
compact. But $M$ is and the limit process as in Theorem \ref{thm of multi-valued
convexity} converges to a point in $V$ since $\nu(\wt p_n)\equiv\vp$ for
any $n$.} there must be a point $\wt q\in\wt V$ and an attained value $\vp+r\vt^*=\wt\mu(\wt q)$ that is not in $\left((\vp+\R\vt^*)\cap\wt\mu(\wt v)\right)^\circ$. Hence there must be a plane, defined by a linear equation
$$
\langle(\xi_0,\vt_0),-\rangle=\langle(\xi_0,\vt_0),\wt\mu(\wt q)\rangle \quad \vt_0\in\toro, \xi_0\in\s,
$$
that bounds on one side the image of the local momentum map at $q$. Observe that  this plane must be transverse to the line $\vp+\R\vt^*$, i.e. $\vt^*(\vt_0)\neq 0$. In particular $\vt_0$ must be non zero. 

The Hamiltonian associated to $(\xi_0,\vt_0)$ has a local maximum at $\wt q$, and hence, by Lemma \ref{4.6} must be exact in $M$. This is inconsistent with the hypotheses of the proposition, given that $\vt_0$ is not zero.

\end{proof}

For simplicity we stated the above result assuming that the 1-forms $\tau_\vt$ are not exact for any $\vt\in\toro$. It is left to the reader to see that this condition can be dropped.

\begin{thm}
Let $M$ be a closed connected symplectic manifold, $G\times_F
T$ a compact connected Lie group acting on $M$, $\wt\mu$ the momentum
map defined from a covering of the manifold $M$. There exists a splitting of $T$ in $T_e\times T_c$ such that $\im\wt\mu\cap(\mathfrak{s}_+^*\times\toro^*)=P\times\toro_c^*$, with $P$ a convex polytope in $\mathfrak s_+^*\times\toro_e^*$. 
\end{thm}

\section{Stability}

It is very natural to ask the following question: How does the image of the momentum map change if the symplectic form is changed into another $T$-invariant symplectic form? One would expect some one-valued Hamiltonians to become multi-valued. In other words there is no reason to believe that, if the Hamiltonian associated to an element $\vt\in\toro$ for a given $T$-invariant symplectic 2-form is exact, then the Hamiltonian associated to the same $\vt$ by another $T$-invariant symplectic 2-form will be exact. 
 
\begin{thm}
Let $(M,\Omega)$ be a compact symplectic manifold, $T$ an $n$-dimensional torus acting on $M$. The subtorus $T_e$, defined as the subgroup whose induced action is Hamiltonian, is constant under small equivariant changes of the symplectic 2-form.
\end{thm}

The topology in the space of 2-forms is the compact-open topology on the space of their coefficients.

\begin{proof} Let $\Omega'=\Omega+\ve$ be a small perturbation of the given form $\Omega$ ($\ve$ is a small closed 2-form). The theorem states that the 1-form $\Omega(X_\vt,-)$ is exact if and only if the 1-form $\Omega'(X_\vt,-)$ is exact.

There is one obvious implication in this statement: non exact Hamiltonians remain non exact Hamiltonians. If $\tau_\vt$ is not exact then it must be true that $\int_\gamma \tau_\vt \neq 0$ for some $\gamma$ non trivial loop in the fundamental group of the manifold. Any change in the symplectic form changes smoothly the form $\tau_\vt$, and so changes smoothly the integral that must hence be non zero for small deformations of $\Omega$ (non exactness is an open property).

What is not obvious is local stability of exact Hamiltonians. In Subsection \ref{MB1-forms} we proved that $\tau_\vt$ is exact if and only if it has a minimum. This is equivalent to the existence of a stationary point $p$ for $\tau_\vt$ (i.e. $\tau_\vt(p)=0$) such that the Hessian matrix of $\tau_\vt$ at $p$ is positive definite; this is equivalent to the fact that the real infinitesimal weights of the isotropy representation at $p$, evaluated at $\vt$, give positive numbers or zero. 

Only the signs of the weights are dependent on the symplectic structure, and their dependence is continuous. So, the existence or not existence of stationary points with the above properties is independent of small equivariant changes of the symplectic 2-form; hence exactness is.
\end{proof}

The above theorem is stated for the case of torus action. The equivalent statement when the acting group is a compact connected Lie group is left to the reader. An immediate corollary of the previous theorem is

\begin{cor}
Let $(M,\Omega)$ be a closed symplectic manifold with a Hamiltonian torus action. The action is Hamiltonian for any other $T$-invariant symplectic structure close to $\Omega$. Also, the same action is Hamiltonian for any other $T$-invariant symplectic form that belongs to the same connected component of $\Omega$ in the vector space of 2-forms.
\end{cor}

The only features in the local expression of the momentum map that are influenced by the symplectic form are: the local parametrization of the polytope at the vertex and the verse of the edges coming out of a vertex. Hence, as long as the weights of the isotropy representation don't change in sign, the only features of the image that can change under perturbation of the symplectic form are the length of the edges of the polytope, not the angles at the vertices.

There are conditions on the possible change in lengths of the edges. It is not true that any such change would return a convex polytope. In fact arbitrary changes in lengths would create a local polytope that ``shifts''.

\section*{Acknowledgments}

I would like to point out that, even though this work is independent of that of Y. Benoist (\cite{bns98}), the statements in that article convinced the author to eliminate a natural requirement in Theorem \ref{thm of multi-valued convexity}, namely Poissonian nature of the action, that made the technical lemmas easier but the statement of the theorem less beautiful. 

I would like to thank S.P. Novikov for posing the problem; F. Cardin, R. Goldin, M. Spera and A.Sikora for discussions that resulted very important in various aspects of this work.

\end{document}